\documentclass[12pt]{amsart}

\usepackage{epsfig}
\usepackage{latexsym}
\theoremstyle{plain}

\newtheorem{thm}{Theorem}

\usepackage[all]{xy}

\theoremstyle{remark}
\newtheorem{exm}{Example}
\newtheorem*{rmk}{Remark}

\usepackage{amssymb}
\usepackage{amsthm}
\usepackage{amsmath}
\usepackage{amsxtra}

\DeclareMathOperator{\GL}{GL}

\DeclareMathOperator{\Frob}{Frob}
\DeclareMathOperator{\Gal}{Gal}

\DeclareMathOperator{\Log}{Log}

\newcommand{\F}{\mathbb F}
\newcommand{\Fqbar} {\overline \F_q}
\newcommand{\Gm}{{\mathbb G}_m}
\newcommand{\Ga}{{\mathbb G}_a}
\newcommand{\C}{\mathbb C}

\newcommand{\Z}{\mathbb Z}
\newcommand{\N}{\mathbb N}

\newcommand{\Ccal}{\mathcal{C}}
\newcommand{\Pcal}{\mathcal{P}}
\newcommand{\Zcal}{\mathcal{Z}}

{.tfm}
{.tfm}

\begin{document}

\title{
Counting colorings on varieties
}

\author{ Fernando Rodriguez-Villegas}
\address{Department of Mathematics 
University of Texas at Austin, TX 78712}

\thanks{ I would like to thank the organizers of the Primeras Jornadas
  de Teor\'{\i}a de N\'umeros in Vilanova i la Geltr\'u for the
  invitation to speak and everyone involved for their generosity and
  hospitality.  Support for this work was supported in part by a grant
  of the NSF }

\email {villegas\@math.utexas.edu}

\maketitle

\section{{}}
\subsection{Introduction}

The goal of this note is to present a combinatorial mechanism for
counting certain objects associated to a variety $X$ defined over a finite
field. The basic example, discussed in \S \ref{2nd-form}, is that of
counting conjugacy classes in $\GL_n(\F_q)$, where $X=\Gm$ (the
multiplicative group). 

We give four different forms of the main formula (which is somewhat
reminiscent of Polya's theory of counting). The principle that emerges
is that in a given setup the counting generating functions for
$X=\bullet$ (a point), $X=\Ga$ (the additive group) and $X=\Gm$ are
related to one another in a simple way. Often one of the cases will be
significantly easier to compute than the others yielding a closed
formula for all three generating functions.  For example, in \S
\ref{applications} we describe how one can go from counting all
matrices in $M_n(\F_q)$, corresponding to $X=\Ga$, to counting
unipotent matrices in $M_n(\F_q)$, corresponding to $X=\bullet$.

None of the special cases considered here are really new; the point is,
instead, to stress the main combinatorial principle. For more general
applications (to quiver and character varieties) we refer the reader
to \cite{HRV} and \cite{HLRV}.

\subsection{The Zeta Function of $X$}

Let $\F_q$ be a finite field with $q$ elements. Fix an algebraic
closure $\Fqbar$ of $\F_q$. For each $r\in \N$ let $\F_{q^r}$ be the
unique subfield of $\Fqbar$ of cardinality $q^r$. Let $\Frob_q\in
\Gal(\Fqbar/\F_q)$ be the Frobenius automorphism $x\mapsto x^q$.  Then
$\F_{q^r}$ is the fixed field of $\Frob_q^r$. 

Let $X$ be an algebraic variety defined over $\F_q$. For each $r \in
\N$ let $N_r(X):=\#X(\F_{q^r})$.  The {\it zeta function} of $X$ is
defined as
\begin{equation}
  \label{zeta}
  Z(X,T):= \exp\left(\sum_{r\geq 1} N_r(X) \frac{T^r}r\right).
\end{equation}
Let $\tilde N_d(X)$ be the number of Frobenius orbits in $X(\Fqbar)$
of size $d$. Then
\begin{equation}
  \label{Nr-tilde}
  N_r(X)= \sum_{d\mid r} d \,\tilde N_d(X).
\end{equation}
We can write the zeta funcion as an Euler product
\begin{equation}
  \label{zeta-euler-prod}
  Z(X,T)=\prod_{d\geq 1}(1-T^d)^{-\tilde N_d(X)}.
\end{equation}

\subsection{Colorings on $X$} 

We consider the following general setup. Let $\Ccal$ be a set, whose
members we  call {\it colors}, and
\begin{equation}
  \label{degree-function}
  |\cdot|: \Ccal \longrightarrow \Z_{\geq 0}
\end{equation}
a function called {\it degree} such that
\begin{enumerate}
\item there are finitely many colors of a given degree;
\item there is a unique color $0\in \Ccal$ of degree $0$.
\end{enumerate}

 A {\it coloring} on $X$ is a map
 \begin{equation}
   \Lambda: X(\Fqbar) \longrightarrow \Ccal.
 \end{equation}
The {\it degree} of $\Lambda$ is defined as
\begin{equation}
  \label{coloring-degree}
  |\Lambda|:= \sum_{x\in X(\Fqbar)} |\Lambda(x)|.
\end{equation}
We will only consider colorings of finite degree, so that
$\Lambda(x)=0$ for all but finitely many $x$.  We let the Frobenius
automorphism act on colorings via
\begin{equation}
  \label{coloring-frob-action}
  \Lambda^{\Frob_q}(x):=\Lambda(\Frob_q(x))
\end{equation}
and say $\Lambda$ is defined over $\F_{q^r}$  if $\Lambda$ is fixed by
$\Frob_q^r$. In this case we will write: $\Lambda$ is a coloring of
$X/\F_{q^r}$. 

Given a pair $(d,\lambda)$, with $d\in \N$ and $\lambda \in \Ccal$ a
non-zero color, we define its {\it multiplicity} $m_{d,\lambda}$ in a
coloring $\Lambda$ of $X/\F_q$ to be the number of Frobenius orbits
$\{x\}$ in $X(\Fqbar)$ of degree $d$ with $\Lambda(x)=\lambda$. Note
that 
\begin{equation}
  \label{degree-mult}
  |\Lambda|=\sum_{d\geq 1,\lambda \neq 0} m_{d,\lambda}\,d|\lambda|.
\end{equation}
We call the combinatorial data $\{m_{d,\lambda}\}$ of multiplicities
the {\it type} of $\Lambda$ and denote it $\tau(\Lambda)$.

\begin{exm}
   Let $\Ccal=\Z_{\geq 0}$ with degree function $|n|=n$. Then a
coloring is an effective $0$-cycle on $X$. The actions of Frobenius
are compatible hence $\Lambda$ is defined over $\F_{q^r}$ if and only
if the corresponding $0$-cycle is.
\end{exm}
\begin{exm}
\label{conj-classes}
If $X=\Gm$ and $\Ccal=\Pcal$ is the set of all partitions of
non-negative integers with $|\lambda|=\lambda_1+\lambda_2+\cdots$ if
$\lambda=\lambda_1\geq\lambda_2\geq\ldots$ then colorings of degree
$n$ are in one-to-one correspondence with conjugacy classes in
$\GL_n(\Fqbar)$ by the Jordan decomposition theorem.  Indeed, to each
coloring $\Lambda$ we may associate the direct sum of Jordan blocks
with eigenvalue $x\in \Fqbar^\star$ and size $\lambda_i$, where
$\Lambda(x)=\lambda_1\geq\lambda_2\geq \ldots$. This correspondence
preserves the action of  Frobenius and therefore colorings defined
over $\F_{q^r}$ are in bijection to conjugacy classes of
$\GL_n(\F_{q^r})$. 

 Similar statements hold for $X=\Ga$ with colorings corresponding to
 conjugacy classes in $M_n(\Fqbar)$ instead.
\end{exm}

We need one more ingredient. Let
$R=\Z[[t_1,\ldots,t_N]][t_1^{-1},\ldots,t_N^{-1}]$ be the ring of
Laurent series with integer coefficients in the variables
$t_1,\ldots,t_N$. We let
\begin{equation}
  \label{weight}
  W:\Ccal \longrightarrow R
\end{equation}
be a function called {\it weight} satisfying $W(0)=1$. 

We define the weight of a coloring $\Lambda$ on $X/\F_q$ as
\begin{equation}
  \label{weight-coloring}
  W(\Lambda):=\prod_{\{x\}}W(\Lambda(x))(t^{d(x)}),
\end{equation}
where $\{x\}$ runs through the Frobenius orbits in $X(\Fqbar)$,
$d(x)=\#\{x\}$ is the degree of $x$ (the size of its Frobenius orbit)
and $t^d:=t_1^d\cdots t_N^d$. Note that $W(\Lambda)$ only depends on
the type $\tau(\Lambda)$:
\begin{equation}
  \label{weight-mult}
  W(\Lambda)=\prod_{d\geq 1,\lambda\neq 0}W(\lambda)(t^d)^{m_{d,\lambda}}.
\end{equation}

We say that $W$ is {\it homogeneous} if for each $\lambda\in \Ccal$ we
have that $W(\lambda) \in R$ is homogeneous of degree $|\lambda|$. In
this case $W(\Lambda)$ is also
homogeneous of degree $|\Lambda|$.

\subsection{Coloring Zeta Function  of $X$} 
Given the coloring data $C=(\Ccal,|\cdot|,W(\cdot))$ we define the
{\it coloring zeta function} of $X$ as the formal power series in
$R[[T]]$
\begin{equation}
  \label{coloring-zeta}
  Z_C(X,t,T):=\sum_{\Lambda} W(\Lambda)T^{|\Lambda|},
\end{equation}
where the sum runs over all colorings of $X/\F_q$.  If $X=\bullet$ (a
point) then the coloring zeta function simply reduces to
\begin{equation}
  \label{coloring-zeta-point}
  Z_C(\bullet,t,T):=\sum_{\lambda \in \Ccal} W(\lambda)T^{|\lambda|}.
\end{equation}

\begin{exm}
 In the {\it standard setup} $C=(\Z_{\geq 0}, |\cdot|,1)$,
  $Z_C(X,T)$ is just the usual zeta function $Z(X,T)$. In particular,
  if $X=\bullet$ then
$$
Z_C(\bullet,T)=\sum_{n\geq 0} T^n= (1-T)^{-1}.
$$
\end{exm}
\begin{exm}
In the {\it partition setup} $C=(\Pcal,|\cdot|,1)$ and if
  $X=\bullet$ then
$$
Z_C(\bullet,T)=\sum_{\lambda\in \Pcal}T^{|\lambda|}=\prod_{d\geq
  1}(1-T^d)^{-1}. 
$$
\end{exm}

\section{{}}

\subsection{First Form}
This form of the main formula is similar to the Euler product
\eqref{zeta-euler-prod} for the usual zeta function (to which it
reduces to in the standard setup).
\begin{thm}
The following identity of generating functions holds
\begin{equation}
  \label{first-form}
  Z_C(X,t,T)=\prod_{d\geq1}Z_C(\bullet,t^d,T^d)^{\tilde N_d}  
\end{equation}
\end{thm}

 \begin{proof}
   Write $Z_C(\bullet,t,T)=1+z(T)$. For $N\in \N$ we have
$$
Z_C(\bullet,t,T)^N=1+\sum_{m\geq 1} N(N-1)\cdots
(N-m+1)\frac{z(T)^m}{m!}
$$
by the binomial theorem. On the other hand by the multinomial theorem
$$
\frac{z(T)^m}{m!}= \sum_{m_\lambda} \prod_{\lambda\neq 0}
\frac{W(\lambda)^{m_\lambda}}{m_\lambda!} T^{m_\lambda|\lambda|}
$$
summed over all sequences of non-negative integers $m_\lambda$ with
$\sum_{\lambda\neq 0} m_\lambda=m$. Putting these two identities
together we get that the coefficient of $T^n$ on the right hand side
of \eqref{first-form} equals
$$
\sum_{m_{d,\lambda}}\;\prod_{d\geq 1,\lambda\neq 0}
\tilde N_d(\tilde N_d-1)\cdots (\tilde N_d-m_d+1)
\frac{W(\lambda)(t^d)^{m_{d,\lambda}}}{m_{d,\lambda}!}
$$
summed over all $m_{d,\lambda}$ sequences of non-zero integers 
satisfying 
$$
n=\sum_{d\geq 1,\lambda \neq 0} m_{d,\lambda}\,d|\lambda|,
$$
where $m_d:=\sum_{\lambda\neq 0} m_{d,\lambda}$.

On the other hand to give a coloring of $X/\F_q$ with multiplicites
$m_{d,\lambda}$ we need to pick $m_d=\sum_{\lambda\neq 0}
m_{d,\lambda}$ Frobenius orbits of size $d$ and color $m_{d,\lambda}$
of them with color $\lambda\neq 0$. There are $\binom{\tilde N_d}{m_d}$
ways of picking the orbits and $m_d!/\prod_{\lambda\neq 0}
  m_{d,\lambda}!$ ways to color them in this way and the weight of
$\Lambda$ is $W(\Lambda)=\prod_{d\geq 1,\lambda\neq
  0}W(\lambda)(t^d)^{m_{d,\lambda}}$. It follows that the
coefficients of $T^n$ on both sides of \eqref{first-form} agree.
 \end{proof}

\subsection{Second Form}
\label{2nd-form}
\begin{thm}
The following identity of generating functions holds
\begin{equation}
  \label{second-form}
  Z_C(X,t,T)=\prod_{m\in \Z^N, d\geq 1}Z(X,t^mT^d)^{v_{d,m}}
\end{equation}
where the exponents $v_{d,m}$ are defined by  the formal identity
\begin{equation}
  \label{Z-prod-exp}
    Z_C(\bullet,t,T)=\prod_{m\in \Z^N, d\geq 1}(1-t^mT^d)^{-v_{d,m}}.
\end{equation}  
\end{thm}
\begin{proof}
  Taking logarithms of both sides of \eqref{first-form} we get
$$
\log Z_C(X,t,T)=\sum_{d\geq 1} \tilde N_d \log Z_C(\bullet,t^d,T^d).
$$
By M\"obius inversion of \eqref{Nr-tilde}
$$
\tilde N_d = \frac 1d \sum_{e\mid d}\mu(e)N_{d/e}.
$$
Taking logarithm of both sides of \eqref{Z-prod-exp} we get
$$
\log Z_C(\bullet,t,T)=-\sum_{m\in \Z^N,d\geq 1} v_{d,m}
\log(1-t^mT^d). 
$$
On the other hand,
$$
T=-\sum_{r\geq 1} \mu(r) \log(1-t^r)
$$
and hence
\begin{eqnarray*}
  \log Z_C(X,t,T) &=& -\sum_{r,s\geq 1} \frac 1{rs} \mu(r) N_s
  \sum_{m\in \Z^N, d\geq 1} v_{d,m}\log(1-t^{msr}T^{dsr})\\
 &=& \sum_{s\geq 1} \frac1s N_s \sum_{m\in \Z^N, \, d\geq 1}
v_{d,m}  t^{sm}T^{sd}\\
&=& \sum_{m\in \Z^N,\, d\geq 1} v_{d,m}\log Z(X,t^mT^d)
\end{eqnarray*}
proving our claim.
\end{proof}
\begin{rmk}
  It is easy to see by induction that the $v_{d,m}$ in
  \eqref{Z-prod-exp} are integers  uniquely determined by
  $Z_C(\bullet,t,T)$. 
\end{rmk}

\begin{exm}
\label{conj-classes-1}
In the {\it partition setup} $C=(\Pcal,|\cdot|,1)$ with
  $X=\Gm$ we have
$$
Z_C(\Gm,T)=\sum_{n\geq 0} C_nT^n,
$$
where $C_n$ is the number of conjugacy classes in $\GL_n(\F_q)$ (see
example \ref{conj-classes}); by \eqref{second-form} this equals
\begin{equation}
\label{conj-class-GL_n}
\prod_{n\geq 1} \left(\frac{1-T^n}{1-qT^n}\right)
\end{equation}
as $Z(\Gm,T)=(1-T)/(1-qT)$. See  \cite{FF}, \cite{G}, \cite{Mc1}.

Similarly, if we take $X=\Ga$ then we get the expression
\begin{equation}
\label{conj-class-M_n}
\prod_{n\geq 1} (1-qT^n)^{-1}
\end{equation}
for the generating funtion for the conjugacy classes in $M_n(\F_q)$
instead. 
\end{exm}

\subsection{Third Form}
This form of the expression for the coloring zeta function is a simple
variant of the second form \eqref{second-form} but it is convenient to
state it separately.

Let
  \begin{equation}
    \label{zeta-u}
    Z_C(u,t,T):=\prod_{m\in \Z^N,d\geq 1} (1-u\,t^mT^d)^{-v_{d,m}}
  \end{equation}
where $u$ is another formal variable and $v_{d,m}$ is as in
\eqref{Z-prod-exp}.   

\begin{thm}
Let
\begin{equation}
  \label{zeta-fact}
Z(X,T)=\prod_i(1-x_iT)^{-n_i},  
\end{equation}
for some $x_i \in \C$ and  $n_i\in \Z$. 
Then with the above notation we have
\begin{equation}
  \label{third-form}
  Z_C(X,t,T)=\prod_iZ_C(x_i,t,T)^{n_i}.
\end{equation}
\end{thm}

\begin{exm}
  In the standard setup  $Z_C(\bullet,T)=(1-T)^{-1}$ so that
  $Z_C(u,T)=(1-uT)^{-1}$ and \eqref{third-form} is simply a
  restatement of \eqref{zeta-fact}.
\end{exm}
\begin{rmk}
  It is known by the work of Dwork that $Z(X,T)$ is a rational
  function of $T$ of the form \eqref{zeta-fact}.
\end{rmk}

\subsection{Fourth Form}

Recall that $R =\Z[[t_1,\ldots,t_N]][t_1^{-1},\ldots,t_N^{-1}]$ is the
ring  of Laurent series in variables $t_1,\ldots,t_N$ with integer
coefficients. Given $Z\in 1+TR[[T]]$ we define, following Getzler
\cite{Ge} 
\begin{equation}
  \label{L}
  \Log(Z):=\sum_{d\geq 1,\, m\in \Z^N}v_{d,m}t^mT^d \in R[[T]],
\end{equation}
where 
$$
    Z=\prod_{m\in \Z^N,\, d\geq 1}(1-t^mT^d)^{-v_{d,m}},
$$
as in \eqref{Z-prod-exp}.

In this section we assume that $X$ is a polynomial count variety; i.e.
$$
N_r(X)=N_X(q^r), \qquad r \in \N,
$$
for some fixed polynomial $N_X\in \Z[q]$. We also assume that one of
the variables in $R$ is $q$. To simplify the notation we relabel the
variables as $q,t_1, \ldots,t_N$ and the exponents as $i\in \Z$ for
$q$ and $m\in \Z^N$ for $t_1,\ldots,t_N$. For example, with this
relabeling   \eqref{zeta-u} becomes
\begin{equation}
    Z_C(u,t,T)=\prod_{i\in \Z,m\in \Z^N,d\geq 1}
    (1-u\,q^it^mT^d)^{-v_{d,i,m}}. 
\end{equation}  

\begin{thm}
  The following identity holds
  \begin{equation}
    \label{fourth-form}
    \Log(Z_C(X,t,T))=N_X(q)\Log(Z_C(\bullet,t,T)).
  \end{equation}
\end{thm}
\begin{proof}
  The claim is a simple consequence of the third form
  \eqref{third-form} of our main formula. If $N_X(q)=\sum_jn_jq^j$
  then
$$
Z(X,T)=\prod_j(1-q^jT)^{-n_j}.
$$
Hence by \eqref{third-form}
\begin{eqnarray*}
  Z_C(X,t,T) &=& \prod_{j\in \Z}Z_C(q^j,t,T)^{n_j}\\
  &=& \prod_{i,j \in \Z, m\in \Z^N, d\geq 1}
  (1-q^{i+j}t^mT^d)^{-n_jv_{d,i,m}}\\
  &=& \prod_{k\in \Z,m\in \Z^N,d\geq
  1}(1-q^kt^mT^d)^{-\sum_{i+j=k}n_jv_{d,i,m}}. 
\end{eqnarray*}
Hence
$$
\Log(Z_C(X,t,T))=\sum_{k\in \Z,m\in \Z^N,d\geq  1}
\sum_{i+j=k}n_jv_{d,i,m}q^kt^mT^d
$$
which equals the right hand side of \eqref{fourth-form}
\end{proof}

\section{{}}
\label{applications}
\subsection{Unipotent matrices}
We consider the coloring data $\Ccal=\Pcal$ with the usual degree
function $|\cdot|$ but with a non-trivial weight function. For all
results and concepts related to partitions our reference will be \cite{Mc},
whose notation we will follow.

 Let $X=\Ga$ and $\Lambda$ a coloring of
$X/\F_q$ corresponding to a conjugacy class $c$ in $M_n(\F_q)$. The
centralizer $z_c$ of $c$ in $G_n:=\GL_n(\F_q)$ has order \cite{Mc}
$$
\prod_{\{x\}}a_{\Lambda(x)}(q^{d(x)}),
$$
where, as before, $\{x\}$ runs through the Frobenius orbits, $d(x)$ is
the size of $\{x\}$ and where for $\lambda\in \Pcal$
\begin{equation}
  \label{a-pol}
  a_\lambda(q):=q^{|\lambda|+2n(\lambda)}b_\lambda(q^{-1}),
\end{equation}
with
\begin{eqnarray}
  \label{a-pol-1}
  n(\lambda)&:=&\sum_{i\geq 1} (i-1)\lambda_i\\
  b_\lambda(q) &:=& \prod_{i\geq 1} \phi_{m_i(\lambda)}(q)\\
  \phi_m(q)&:=&(1-q)(1-q^2)\cdots (1-q^m)
\end{eqnarray}
and, finally, $m_i(\lambda)$ is the multiplicity of $i$ in
$\lambda$. 

It follows that if we define our weight function as
$$W
(\lambda):=a_\lambda(q)^{-1}\in R
$$
 where $R=\Z[[q]][q^{-1}]$ then 
\begin{equation}
  \label{centralizer}
  W(\Lambda)=\frac 1 {|z_c|}.
\end{equation}
Consequently, if we now take $X=\Gm$ then
$$
  \sum_{\Lambda/\F_q,|\Lambda|=n} W(\Lambda)=1
$$
and therefore
\begin{equation}
  \label{centr-zeta-Gm}
  Z_C(\Gm,q,T)=\sum_{n\geq 0}T^n=(1-T)^{-1}.
\end{equation}

Applying \eqref{fourth-form} to this situation we find that
\begin{equation}
  \label{centr-zeta-point}
  Z_C(\bullet,q,T)=\prod_{n\geq 0}(1-q^nT),
\end{equation}
since $N_{\Gm}=q-1$ and $(q-1)^{-1}=-(1+q+q^2+\cdots)$; sincen
 $N_{\Ga}=q$ by  \eqref{fourth-form} we also have the identity
 \begin{equation}
   \label{centr-zeta-Ga}
   Z_C(\Ga,q,T)=\prod_{n\geq 1}(1-q^nT).
 \end{equation}

On the other hand, for $X=\Ga$ we have
$$
\sum_{|\lambda|=n}\frac 1{a_\lambda(q)}= \sum_{\Lambda/\F_q,|\Lambda|=n}
 W(\Lambda)=\frac{|M_n(\F_q)|}{|G_n|}
 =\frac{q^{\tfrac12n(n+1)}}{(q^n-1)(q^{n-1}-1)\cdots (q-1)}
$$
and we have therefore proved the following identity of Euler
\begin{equation}
  \label{euler}
\sum_{n\geq 0} \frac{q^{\tfrac12n(n+1)}T^n}{(q^n-1)(q^{n-1}-1)\cdots
  (q-1)}=\prod_{n\geq 1}(1-q^nT).  
\end{equation}

If, instead,  $X=\bullet$ we obtain
  $$
 \sum_{\Lambda/\F_q,|\Lambda|=n}
 W(\Lambda)=\frac{u_n}{|G_n|},
 $$
 where $u_n$ is the number of unipotent matrices in $G_n$.
 Combining \eqref{centr-zeta-point} with Euler's identity
 \eqref{euler} with $T$ replaced by $T/q$  we find
 $$
 \frac{u_n}{|G_n|}=
 \frac{q^{\tfrac12n(n+1)-n}}{(q^n-1)(q^{n-1}-1)\cdots (q-1)}
$$
we deduce the known result $u_n=q^{n^2-n}$ (see \cite{St} for a general
result on the number of unipotent elements in linear algebraic groups
over finite fields). 

\subsection{Commuting pairs of matrices}
We now consider a weight function arising from the centralizer algebra
$\Zcal_A$ of a matrix $A\in M_n(\F_q)$. It is known that 
$$
\dim_{\F_q}(\Zcal_A)=\sum_{\{x\}}d(x)\langle \Lambda(x),\Lambda(x)\rangle
$$
where for a partition $\lambda$ we define $\langle \lambda, \lambda
\rangle:= |\lambda| + 2n(\lambda)$. 

Since $|\Zcal_A|$ only depends on the conjugacy class $[A]$ of $A$, we
can count ordered pairs of commuting matrices in $M_n(\F_q)$ as follows
$$
\gamma_n:=\#\{A,B \in M_n(\F_q) \;\mid \; AB=BA\}=
\sum_{[A]} \#[A] \,|\Zcal_A|, 
$$
where $[A]$ runs through the conjugacy classes in $M_n(\F_q)$.
Hence if we define
$$
W(\lambda):=\frac{q^{\langle\lambda,\lambda\rangle}}{a_\lambda(q)}
$$
then 
$$
W(\Lambda)=\frac{\#[A]}{|G_n|}|\Zcal_A|,
$$
where $[A]$ corresponds to the coloring $\Lambda$  on $\Ga/\F_q$,
and therefore
$$
\gamma_n=|G_n|\sum_\Lambda  W(\Lambda).
$$
 Consequently,
\begin{equation}
  \label{comm-zeta}
Z_C(\Ga,q,T)= \sum_{n\geq 0} \frac{\gamma_n}{|G_n|}T^n.
\end{equation}

On the other hand,
\begin{eqnarray*}
  Z_C(\bullet,q,T) &=& \sum_\lambda
  \frac{q^{\langle\lambda,\lambda\rangle}}{a_\lambda(q)}
  T^{|\lambda|}\\
  &=& \sum_\lambda \frac{T^{|\lambda|}}{b_\lambda(q^{-1})}\\
  &=& \prod_{i\geq 1}\sum_{m_i\geq 0}
  \frac{T^{im_i}}{\phi_{m_i}(q^{-1})}\\
    &=& \prod_{i,n\geq 1}(1-q^nT^i)
\end{eqnarray*}
using Euler's identity \eqref{euler}. 
 Applying \eqref{third-form} we recover (in an equivalent form)  a
 result of Fine and Feit \cite{FF}
 \begin{equation}
   \label{fine-fite}
\sum_{n\geq 0} \frac{\gamma_n}{|G_n|}\;T^n=\prod_{i,n\geq
  1}(1-q^{n+1}T^i).
 \end{equation}

Similarly, we obtain
$$
\sum_{n\geq 0} \frac{\gamma'_n}{|G_n|}\;T^n= Z_C(\Gm,q,T),
$$
where
$$
\gamma'_n:=\#\{A \in \GL_n(\F_q), B \in M_n(\F_q) \;\mid \; AB=BA\}.
$$
Again by \eqref{third-form} we find
$$
Z_C(\Gm,q,T)= \prod_{i,n\geq 1}
\left(\frac{1-q^{n+1}T^i}{1-q^nT^i}\right) = \prod_{i\geq 1}
(1-qT^i)^{-1}.
$$
We now recognize this generating series as \eqref{conj-class-M_n}
and conclude that $\gamma_n'/|G_n|$ is the number of conjugacy classes
in $M_n(\F_q)$. This, in fact, can be proved directly by a simple
application of Burnside's lemma to $\GL_n(\F_q)$ acting on $M_n(\F_q)$
by conjugation. By our main combinatorial principle, this means that
we can run the argument backwards and prove \eqref{fine-fite}
starting from \eqref{conj-class-M_n}.


\begin{thebibliography}{1} 
  
\bibitem{HRV} T. Hausel and F. Rodriguez-Villegas, {\it On the
    E-polynomial of certain     character varieties}, in preparation. 

\bibitem{HLRV} T. Hausel, E. Letellier and F. Rodriguez-Villegas, in
  preparation. 

\bibitem{FF}  W. Feit and N. J. Fine, {\it  Pairs of commuting
    matrices over a finite field},  Duke Math. J.  {\bf 27}  (1960)
  91--94. 
\bibitem{Ge}
E. Getzler,
{\it Resolving mixed Hodge modules on configuration spaces},
Duke Math. J. {\bf 96} (1999), 175--203

\bibitem{G} J. A. Green,  {\it The characters of the finite general
    linear groups},  Trans. Amer. Math. Soc.  {\bf 80}  (1955),
  402--447. 
\bibitem{Mc1} I. G. Macdonald, 
{\it Numbers of conjugacy classes in some finite classical groups}, 
Bull. Austral. Math. Soc. {\bf 23} (1981), 23--48.

\bibitem{Mc}  I. G. Macdonald, {\it Symmetric functions and Hall
    polynomials}, Second edition. With contributions by
    A. Zelevinsky. Oxford Mathematical Monographs. Oxford Science
    Publications. The Clarendon Press, Oxford University Press, New
    York, 1995.

\bibitem{St} 
Robert, Steinberg, {\it 
Endomorphisms of linear algebraic groups} 
Memoirs of the AMS, No. {\bf 80}
AMS, Providence, R.I. 1968
\end{thebibliography}
 \end{document}